\documentclass[12pt]{article}
\usepackage{amsmath, stackrel}
\usepackage{amssymb}
\usepackage{dsfont}
\usepackage{cite}
\newsymbol \blackbox 1004
\newcommand{\eh}{\hfill}\newlength{\sperr}

\newenvironment{proof}{{\settowidth{\sperr}{\bf\rm
Proof}%
\par\addvspace{0.3cm}\noindent\parbox[t]{1.3\sperr}
{\bf\rm P\eh r\eh o\eh o\eh f\eh }%
}}{\nopagebreak\mbox{}
$\blackbox$\par\addvspace{0.3cm}}

\def\nn{\nonumber}
\def\ov{\overline}

\def\b{\beta}

\def\s{\sigma}
\def\la{\lambda}

\def\wt{\widetilde}
\def\BN{{\mathbb N}}
\def\BC{{\mathbb C}} 
\def\BR{{\mathbb R}} 
\def\clb{{\mathcal B}}

\def\clh{{\mathcal H}}
\def\cli{{\mathcal I}}
\def\clk{{\mathcal K}}
\def\clw{{\mathcal W}}

\def\mfa{{\mathfrak A}}

\newtheorem{Pa}{Paper}[section]
\newtheorem{Tm}[Pa]{{\bf Theorem}}
\newtheorem{La}[Pa]{{\bf Lemma}}
\newtheorem{Cy}[Pa]{{\bf Corollary}}
\newtheorem{Rk}[Pa]{{\bf Remark}}

\newtheorem{Ee}[Pa]{{\bf Example}}

\newtheorem{Pn}[Pa]{{\bf Proposition}}

\newcommand{\E}{\mathrm{e}}
\newcommand{\I}{\mathrm{i}}

\author{Alexander Sakhnovich \footnote{The author was supported by the Austrian Science Fund (FWF) grant  DOI: 10.55776/Y963}}
\title{GBDT, multiplicative integrals and linear similarity}

\date{}
\begin{document}

\maketitle


\begin{flushleft}
Faculty of Mathematics,
University
of
Vienna, \\
Oskar-Morgenstern-Platz 1, A-1090 Vienna,
Austria\\
E-mail: oleksandr.sakhnovych@univie.ac.at
\end{flushleft}

\begin{abstract} A GBDT version of the  B\"acklund-Darboux transformation
for a non-isospectral canonical system is considered. 
Applications to  multiplicative integrals and their limit values, to characteristic matrix functions and 
to linear similarity problems are obtained. Some interesting examples are constructed as well.
\end{abstract}

{\it 2020 Mathematics Subject Classification.} 34A05, 34L40, 45H05, 47A48.

{\it Keywords.} {Non-isospectral canonical system, B\"acklund-Darboux transformation, GBDT, multiplicative integral, 
characteristic matrix function, linear similarity}.

\section{Introduction}\label{sec0}
\setcounter{equation}{0}
The classical canonical system is given by the equation
\begin{equation} \label{1.1-}
w_x(x,\la)=\I \la J H(x) w(x, \la) ,
\end{equation}
where $w_x= \frac{d}{d x}w$, $\I$ stands for the imaginary unit ($\I^2=-1$), $J$ and  $H(x)$ are $m \times m$
matrices:
\begin{equation} \label{1.2}
H(x) \geq 0, \quad J=J^*=J^{-1},
\end{equation}
and the independent variable $\la$ is a so called spectral parameter.
Here, $m\in \BN$ (where $\BN$ is the set of the positive integer values), $J^*$ stands for 
the  conjugate transpose for $J$ and $H(x)\geq 0$ means that $H(x)$  is Hermitian and the eigenvalues of the Hamiltonians $H(x)$
are nonnegative at almost all $x$.

Canonical systems have been actively studied in analysis (see, e.g., \cite{dBr, GoKr, KLang, Lang, Su} 
and the references therein) and some
basic results for the case $m=2$ were obtained by M. G. Krein and L. de Branges.   The case $m >2$ 
have been actively studied as well (see \cite{ArD, dB2, SaA3, ALSinvpr, ALSarx, SaSaR, SaL3} and the references
therein) but this case is more complicated. In fact, we have 
\begin{equation} \label{1.2+}
m=2p, \quad J=\begin{bmatrix} 0 & I_p\\I_p &0\end{bmatrix}
\end{equation}
for the classical canonical systems,
where $I_p$ is the $p\times p$ identity matrix.

When $\la$ depends on $x$, we speak about a non-isospectral canonical system. In our case, we have $\la=(z-x)^{-1}$,
where $z$ is a ``hidden" independent  of $x$ spectral parameter. Writing $w(x,z)$ instead of $w\big(x, (z-x)^{-1}\big)$,
we rewrite \eqref{1.1-} in the form
\begin{equation} \label{1.1}
w_x(x,z)=\I  (z-x)^{-1}J H(x) w(x, z) .
\end{equation}
Further, we require that \eqref{1.2} holds but \eqref{1.2+} is not necessary. 
When the  interval $\cli \subset \BR$ (where system \eqref{1.1} is considered) is fixed, the 
$m \times m$ fundamental solution
$W(x,z)$ of (\ref{1.1}) is normalised by the condition
\begin{equation} \label{1.2'}
W(\xi,z)=I_m \quad {\mathrm{for\,\, some\,\, fixed}} \quad \xi\in \cli.
\end{equation}
{\it We always assume that the entries of $H(x)$ are locally integrable, that is, integrable on the closed intervals belonging to $\cli$.}

In several important cases, system \eqref{1.1} on the interval $[\xi, \zeta]$, with  $W(\zeta,z)$ analytic for $z\not\in [\xi,\ell]$ ($\ell\leq \zeta$),
generates a Riemann--Hilbert problem
\begin{align}\label{0.1}&
W_+(\zeta,s)=W_-(\zeta,s)R(s)^2 \quad (\xi \leq s \leq \ell), \quad
{\mathrm{where}}
\\ \label{0.2}&
W_+(\zeta,s):= \lim_{\eta \to +0} W(\zeta, s+\I \eta), \quad W_-(\zeta,s):=
\lim_{\eta \to +0} W(\zeta, s-\I \eta).
\end{align}
This problem is of interest in the random matrix theory. In particular,
the Markov parameters appearing in the series representation of $W(\zeta, z)$ are essential for the random matrices
problems \cite{Dei}.  Ernst-type equation with auxiliary linear systems of the \eqref{1.1} form was studied in 
\cite{ALSgrav}. The interconnections (see   \cite{SaL0, SaL4}) between systems \eqref{1.1}, multiplicative integrals and characteristic
matrix valued functions (further called  characteristic
matrix  functions or simply characteristic
functions) will be used  in the present paper.

Interrelated approaches of B\"acklund-Darboux transformations, dressing, commutation and double commutation methods 
(see, e.g., \cite{C, CI, Dei0, Ge, GeT, MS, Mar, SaA3, 
ALSgrav, SaSaR} and the references therein)
are used in order to construct coefficients and solutions of the families of so called transformed systems in terms of the
simpler coefficients and solutions of an initial system. GBDT version of B\"acklund-Darboux transformation for the
non-isospectral canonical system \eqref{1.1} was constructed in \cite{ALS07}. Here, we apply a simple but useful modification
of the GBDT from \cite{ALS07}.

We use GBDT in order to construct families of more complicated cases  of multiplicative integrals and their limit values and
of the so called triangular models and their characteristic functions given some initial simple case.
The  mentioned above families have (in a certain way) the  properties similar to those of the initial
multiplicative integrals and triangular models, respectively.

The paper consists of 5 sections, the first one being this introduction. In the second (preliminary) section,  we introduce
GBDT for the system \eqref{1.1}. The second section is based mostly on our work  \cite{ALS07}. In Section \ref{sec3}, we study
multiplicative integrals generated by systems \eqref{1.1} and their limit values.  In Section \ref{CharF}, we consider
triangular integral operators (triangular models) and characteristic functions. Section \ref{secEx} is dedicated to explicit examples.


{\bf Notations.} The set of real numbers is denoted by $\BR$, the set of  complex numbers (complex plane) is denoted by $\BC$, $\Re(z)$
stands for the real and $\Im(z)$ for the imaginary part of $z$.
The upper half-plane of  $z\in \BC$, where $\Im (z)>0$,  is denoted by $\BC_+$ and
the lower half-plane of $z \in \BC$, where $\Im (z)<0$, is denoted by $\BC_-$.  The notation  $\ov{g}$ stands for the complex conjugate of $g\in \BC$.
The set of $m\times n$ matrices
with complex valued entries is denoted by $\BC^{m\times n}$ and we put  $\BC^m=\BC^{m\times 1}$. The notation $\s(B)$ stands for the spectrum of the matrix $B$ and
$I$ is the identity operator.
The set of bounded operators acting from some Hilbert space $\clh_1$ into Hilbert space $\clh_2$ is denoted by $\clb(\clh_1,\clh_2)$,
we set $\clb(\clh_1):=\clb(\clh_1,\clh_1)$, and
$\clk^*$ stands for the operator adjoint to the operator $\clk$. Finally, $L^2_k(a,b)$ is the Hilbert space of the square-integrable, $k$-dimensional
vector functions on $(a,b)$.

\section{Preliminaries: GBDT for system \eqref{1.1}}\label{sec2}
\setcounter{equation}{0}
 The methods developed in \cite{SaA2,SaA3} for
non-isospectral problems and canonical systems \eqref{1.1-}, respectively,
have been used in \cite{ALS07} in order to construct GBDT for
the non-isospectral canonical system \eqref{1.1}. The corresponding
procedure is described below.

Let an initial system \eqref{1.1}, such that \eqref{1.2} holds, be given on an interval $\cli\subset \BR$.  Let us fix a point $\xi\in \cli$, an integer $n>0$, and $n \times n$ parameter matrices
$A(\xi)$ and $S(\xi)$. Fix also an $n \times m $ parameter matrix
$\Pi(\xi)$ so that the matrix identity
\begin{equation} \label{1.3}
A(\xi)S(\xi)-S(\xi)A(\xi)^*=\I \Pi(\xi)J \Pi(\xi)^* \quad \big(S(\xi)=S(\xi)^*\big)
\end{equation}
holds. Introduce now matrix  functions  $A(x)$, $S(x)$ and $\Pi(x)$
by their values at $x=\xi$ and equations
\begin{align} & \label{1.4}
A_x=A^2, \quad \Pi_x=- \I A \Pi J H,
\\ & \label{1.5}
S_x= \Pi J H J^* \Pi^*- (AS+SA^*),
\end{align}
where the argument $x$ is omitted  for brevity.
Then, it can be checked by direct differentiation that the matrix
identity
\begin{equation} \label{1.6}
A(x)S(x)-S(x)A(x)^*=\I \Pi(x) J \Pi(x)^*
\end{equation}
holds for each $x$. (We also have $S(x)=S(x)^*$.) Notice that the equation $A_x=A^2$ is motivated by
the similar equation $\la_x= \la^2$ for  the $x$-dependent spectral parameter $\la=(z-x)^{-1}$
because $A$ can be viewed as  a generalised spectral parameter
(compare formulas (3) and (4) for $\la$ and $A$ in \cite{SaA2} ). In the points of invertibility of $S(x)$ we 
introduce a transfer matrix function in the Lev Sakhnovich form
\cite{SaSaR, SaL1, SaL3} (modified for the nonisospectral case \eqref{1.1} with a ``hidden" independent parameter $z$):
\begin{equation} \label{1.7}
w_A(x, z)=I_m-\I J \Pi(x)^*S(x)^{-1}(A(x)- \la(x,z) I_n)^{-1} \Pi(x), \quad \la=(z-x)^{-1}.
\end{equation}
This transfer matrix function has an important $J$-property
\cite{SaL1}:
\begin{equation} \label{1.7'}
w_A(x, \ov{z})^*Jw_A(x, z)=J, \quad \mathrm{i.e.,} \quad
w_A(x,z)^{-1}=Jw_A(x, \ov{z})^*J.
\end{equation}
We set
\begin{equation} \label{1.10}
v(x, z)=w_0(x)^{-1}w_A(x,z).
\end{equation}
where, in  case of the invertible matrices $A(x)$, $w_0(x)$ is defined by the relations
\begin{equation} \label{1.11}
w_0(x):=w_A(x, \infty)=I_m-\I J \Pi(x)^*S(x)^{-1}A(x)^{-1} \Pi(x)
\end{equation}
(see \cite[Remark 1]{ALS07}). Hence, according to \eqref{1.7'} and \eqref{1.10}, $v(x,z)$ is invertible and we
also have
\begin{equation} \label{1.12}
w_0(x)^{-1}=Jw_A(x, \infty)^*J=I_m+\I J \Pi(x)^*\big(A(x)^*\big)^{-1} S(x)^{-1}\Pi(x).
\end{equation}
\begin{Rk}\label{A(x)} If $A(x)$ is well-defined and invertible on some
interval $\cli$, we rewrite $A_x=A^2$ in the form $-A^{-1}A_xA^{-1}=-I_n$ or, equivalently,
in the form $(A^{-1})_x=-I_n$. Thus, $A(x)^{-1}=B-xI_n$ $(B\in \BC^{n\times n})$
on $\cli$ and we have
\begin{equation} \label{1.17}
A(x)=(B-xI_n)^{-1} \quad {\mathrm{for\, all}} \quad x\in \cli, \quad \s(B)\cap \cli=\emptyset.
\end{equation}
It is also immediate that $A_x=A^2$  and  $A(x)$ is well defined and invertible on $\cli$ if \eqref{1.17} holds.
\end{Rk}
In the case \eqref{1.17}, we have
\begin{equation} \label{1.17+}
(A(x)-\la(x,z)I_n)^{-1}=(x-z)(B-xI_n)(B-zI_n)^{-1},
\end{equation}
that is, the resolvent $(A(x)-\la(x,z)I_n)^{-1}$ is well defined for $z\not\in \s(B)$
Now, the main result on GBDT for system \eqref{1.1} (based on \cite[Theorem 1]{ALS07}) takes 
the following form.
\begin{Tm} \label{TmBDT} Suppose that $W(x,z)$ is the fundamental solution
of the initial system \eqref{1.1} on some interval $\cli \subset \BR$ and that relations \eqref{1.2}, \eqref{1.2'}, \eqref{1.3}, \eqref{1.17}, $z, \ov{z}\not\in \s(B)$ and 
\begin{align} & \label{1.14-}
 \Pi_x=- \I A \Pi J H,
\quad
S_x= \Pi J H J^* \Pi^*- (AS+SA^*)
\end{align}
are valid on $\cli$. Assume that $S(x)$ is invertible on $\cli$ and  let $v(x,z)$ be given by \eqref{1.10} and \eqref{1.11}.
Then, the matrix
function
\begin{equation} \label{1.14}
\wt w(x, z)=v(x, z) W(x, z)
\end{equation}
is well defined and satisfies the transformed $($GBDT-transformed$)$ system
\begin{equation} \label{1.13}
\frac{d}{dx} \wt w(x,z)=\I (z-x)^{-1} J \wt H(x) \wt w(x,z),
\end{equation}
where the transformed Hamiltonian $\wt H(x)$ has the form
\begin{equation} \label{1.15}
\wt H(x)=w_0(x)^*H(x)w_0(x).
\end{equation}
Thus,
the fundamental solution $\wt W(x,z)$ of the transformed system \eqref{1.13}
$($normalised by $\wt W(\xi,z)=I_m)$ is given by the
formula
\begin{equation} \label{1.16}
\wt W(x, z)=v(x, z) W(x, z)v(\xi, z)^{-1} \quad (0 \leq x \leq l).
\end{equation}
\end{Tm}
Let us consider a useful condition of invertibility of $S(x)$.
\begin{Pn}\label{PnInv}
Suppose that some initial system \eqref{1.1} satisfying \eqref{1.2} is given on the closed interval $\cli=[a,b]\subset\BR$ and put $\xi=a$ here. Assume that  the triple 
of parameter matrices $ \{A(\xi), S(\xi), \Pi(\xi)\}$ 
$(S(\xi)=S(\xi)^*)$ satisfies the matrix identity \eqref{1.3}, that the relations
\eqref{1.17} and \eqref{1.14-} hold and that $S(\xi)>0$. Then, $S(x)>0$ on $\cli$ and  $\| S(x)^{-1}\|$ is bounded on $\cli$.
If all the conditions of the proposition are valid but $\cli=[a,b)$, we have 
$S(x)>0$ on $\cli$ and  $\| S(x)^{-1}\|$ is bounded on all the closed intervals belonging
to $\cli$.
\end{Pn}
\begin{proof}. We introduce matrix functions $Q(x)$ and $K(x)$ by the relations 
\begin{equation} \label{1.18}
Q(x)=K(x)S(x)K(x)^*, \quad {\mathrm{where}} \quad K_x(x)=K(x)A(x), \,\, K(a)=I_n.
\end{equation}
It follows from \eqref{1.5} and \eqref{1.18} that 
\begin{align}\nn
Q_x(x)&=K(x)\big(S_x(x)+ A(x)S(x)+S(x)A(x)^*)K(x)^*
\\  \label{1.19} &
=K(x)\Pi(x)JH(x)J^*\Pi(x)^*K(x)^*\geq 0, \quad Q(a)=S(a)>0.
\end{align}
It is easy to see that $K(x)$ is invertible. Therefore, \eqref{1.18} and \eqref{1.19} yield
\begin{align}& \label{1.20}
S(x)=K(x)^{-1}Q(x)(K(x)^{-1})^*\geq K(x)^{-1}S(a)(K(x)^{-1})^*>0, 
\\ & \label{1.21}
 S(x)^{-1}\leq K(x)S(a)^{-1}K(x)^*.
\end{align}
If $\cli=[a,b]$, $\|K(x)\|$ is bounded on $[a,b]$ and formula \eqref{1.21} implies that
$\| S(x)^{-1}\|$ is bounded on $[a,b]$ as well.
The statement of the proposition for the case $\cli=[a,b)$ easily follows
from the considerations above.
\end{proof}
\begin{Rk}
Since $A(x)$, $\Pi(x)$ and $S(x)$ are continuous, we have the local boundedness of $w_0(x)$ and the local integrability 
of $\wt H(x)$ under conditions of Theorem \ref{TmBDT}.
According to Proposition \ref{PnInv} $($under stronger conditions above$)$,  we also have the boundedness of the $\|w_0(x)\|$ on $[a,b]$
and one may easily derive certain estimates.
\end{Rk}
\section{GBDT and multiplicative integrals}\label{sec3}
\setcounter{equation}{0}
{\bf 1.} In this section, we assume again that
\begin{align}\label{2.0} &
\cli=[a,b], \quad \xi=a.
\end{align}
Fundamental solutions of systems \eqref{1.1} on $\cli$ admit multiplicative integral representations  (see \cite{SaL4} and the references therein):
\begin{align}\label{2.1} &
W(x,z)=\overset{x}{\overset{\curvearrowleft}{\underset{a}{\int}}}\E^{\frac{\I J}{z-t}d\tau(t)}=\overset{x}{\overset{\curvearrowleft}{\underset{a}{\int}}}\E^{\frac{\I J}{z-t}H(t)dt}
\quad (x\in \cli),
\end{align}
where $\tau(x)=\int_0^x H(t)dt$. We note that \eqref{1.1} implies that $W(x,z)$ is $J$-expanding in the upper half-plane $\Im(z)>0$ and $J$-contractive
in the lower half-plane $\Im(z)<0$, that is
\begin{align}\label{2.2}
W(x,z)^*JW(x,z) \geq J \,\, {\mathrm{for}} \,\, \Im(z)>0, \quad W(x,z)^*JW(x,z) \leq J \,\, {\mathrm{for}} \,\, \Im(z)<0.
\end{align}
The multiplicative representation of the $J$-contractive matrix functions have been studied in the seminal work \cite{Pot}.
\begin{Rk}\label{RkLish}
Formula \eqref{2.1} implies also that $W(x,z)$ is analytic for $z\not\in [a,b]$ and $\lim_{z\to \infty} W(x,z)=I_m$.
If some $m\times m$ matrix function $W(z)$ has these properties, satisfies \eqref{2.2} and $\sup_{z\in (\BC-\BR)} \|W(z)\|<\infty$,
then $W(z)$ admits multiplicative representation \cite[Theorem 4.3]{SaL5} $($see also \cite[Section 3.2]{SaL3}$):$
\begin{align}\label{2.2+} &
W(z)=\overset{b}{\overset{\curvearrowleft}{\underset{a}{\int}}}\E^{\frac{\I J}{z-t}H(t)dt} \qquad (H(t)\geq 0).
\end{align}
\end{Rk}
The limit values $W_{\pm}(x,z)$ of multiplicative integrals (introduced in \eqref{0.2}) have been actively studied.
The corresponding results for multiplicative integrals are considered as certain analogs of the well-known Plemelj
formulas for the usual integrals. In particular, see \cite{SaL0} and the references therein for the important case
$J=I_m$.  Our Theorem \ref{TmBDT}   yields the next proposition.
\begin{Pn}\label{PnMult} Let  $W(x,z)$ admit representation \eqref{2.1}, let relations \eqref{1.2},  \eqref{1.3}, \eqref{1.17},
\eqref{1.14-}, and \eqref{2.0} hold, and  assume that $S(x)$ is invertible on $[a,b]$
Let the limit value $W_+(\ell,s)$ $(W_-(\ell,s))$, where $0<s\leq \ell\leq b$, 
exist. Then, the limit value $\wt W_+(\ell,s)$ $(\wt W_-(\ell,s))$  of the multiplicative integral
\begin{align}\label{2.3} &
\wt W(\ell,z)=\overset{\ell}{\overset{\curvearrowleft}{\underset{a}{\int}}}\E^{\frac{\I J}{z-t}\wt H(t)dt},
\end{align}
where $\wt H$ is given by \eqref{1.15}, exists
and we have 
\begin{align}\label{2.4}
\wt W_+(\ell,s)=v(\ell, s) W_+(\ell, s)v(a, s)^{-1} \quad (\wt W_-(\ell,s)=v(\ell, s) W_-(\ell, s)v(a, s)^{-1}).
\end{align}
\end{Pn}
Here, one  may  use formulas \eqref{1.17} and \eqref{1.17+}  in order to see that
$w_A(\ell,z)$, $v(\ell,z)$ and $v(a,z)^{-1}$ are well defined and continuous in some  neighbourhood of $z=s$.
Therefore, Theorem \ref{TmBDT}  implies, indeed, Proposition \ref{PnMult}.
\begin{Rk} \label{RkF}  The  generalised for the matrix function case Fatou theorem $($see, e.g., \cite[Lemma 1.1]{SaL0}$)$ is essential in the theory
of limit values of $W(z)$. It states that, for the matrix function $W(z)$
analytic $($and with bounded norm$)$ in $\BC_+$, the limit $W_+(s)$ exists almost everywhere on $\BR$.
$($The generalised  Fatou theorem is valid also for the infinite matrix case, where the limit is substituted by
the strong limit.$)$ Clearly, a similar result is valid in the lower half-plane $\BC_-$. 
\end{Rk}
In view of Proposition \ref{PnMult} and Remark \ref{RkF} (and under the conditions of Proposition \ref{PnMult} and Remark \ref{RkF})
the limit values $\wt W_+(\ell,s)$  exist and are given by \eqref{2.4} almost everywhere on $\BR$. A similar result under analogous conditions is valid for $\wt W_-(\ell,s)$.

In a similar to the proof of Proposition \ref{PnMult} way, Theorem \ref{TmBDT}  yields  the boundedness of $\wt W(b,z)$.
\begin{Pn}\label{PnB} Let the norm of $W(b,z)$ of the form \eqref{2.1}   be bounded in $\BC_+$ $(\BC_-)$. 
Let relations \eqref{1.2},  \eqref{1.3}, \eqref{1.17},
\eqref{1.14-}, and \eqref{2.0} hold, and  assume that $S(x)$ is invertible on $[a,b]$.
Then,  the norm of  $\wt W(b,z)$ of the form \eqref{2.3} is bounded in $\BC_+$ $(\BC_-)$ 
as well.
\end{Pn}

{\bf 2.} An important subclass of  the Hamiltionians $H(x)$ of the form
\begin{align}\label{2.5}
H(x)=\b(x)^*\b(x), \quad \b(x)J\b(x)^*\equiv 0 \quad \big(\b(x)\in\BC^{k\times m}\big),
\end{align}
and the corresponding multiplicative integrals \eqref{2.1} have been studied in \cite{SaL4}. In particular,
the following lemma was proved in \cite{SaL4} as Lemma 2.1:
\begin{La}\label{La} Let $H(x)=\b(x)^*\b(x)$ and let the $k \times m$ matrix function $\b(x)$ be continuous on $[a,b]$. Assume that the relation
\begin{align}\label{2.6}
\sup_{a\leq t<x\leq b}\| (x-t)^{-1}\b(x)J\b(t)^*\| \leq M 
\end{align}
holds for some $M>0$. Then, there is a limit
\begin{align}\label{2.7}
V(x,s)=\lim_{\eta \to +0}\big(W(x, s+\I \eta)-W(x, s-\I \eta)\big) \quad (a\leq s\leq x\leq b),
\end{align}
and this limit $V(x,s)$ satisfies $($for some $M_1>0)$ the inequality
\begin{align}\label{2.8}
\sup\|V(x,s)\| \leq M_1 \quad (a\leq s\leq x\leq b).
\end{align}
\end{La}
Conditions of Lemma \ref{La} imply that $\b(x)J\b(x)^*\equiv 0$, that is, the second equality in \eqref{2.5} holds.
In view of \eqref{1.7'} and \eqref{1.11}, we have
\begin{align}\label{2.9}
w_0(x)Jw_0(x)^*=J \quad (x\in [a,b]).
\end{align}
According to \eqref{1.15} and \eqref{2.9},  formula \eqref{2.5} for $H(x)$ yields a similar  to \eqref{2.5}
representation for the transformed Hamiltonian $\wt H(x)$:
\begin{align}\label{2.10}
\wt H(x)=\wt \b(x)^*\wt \b(x), \quad \wt \b(x)J\wt \b(x)^*\equiv 0, \quad {\mathrm{where}} \,\, \wt\b(x):=\b(x)w_0(x).
\end{align}
We will need relations \cite[(2.8) and (2.9)]{ALS07} :
\begin{align} \label{2.11}&
\frac{d}{dx}w_0(x)=  G_0(x)  w_0(x),\\
\label{2.12}&
 G_0=-J (\I \Pi^*S^{-1} \Pi- H J \Pi^* S^{-1} \Pi+
\Pi^* S^{-1} \Pi J^* H).
\end{align}
Let us rewrite \eqref{2.11} in the form
\begin{align} \label{2.13}&
w_0(x)=w_0(t)+\int_t^x  G_0(t)  w_0(t)dt \quad (t<x).
\end{align}
If \eqref{2.5} holds and $\b(x)$ is continuous, then $G_0(x)$ given by \eqref{2.12} (under conditions
of Theorem \ref{TmBDT}) is continuous on $[a,b]$ as well. Since $w_0(x)$ and $G_0(x)$ are continuous on $[a,b]$,
the norm of $G_0(x)w_0(x)$ is bounded on $[a,b]$. Hence, formula \eqref{2.13} implies that
\begin{align} \label{2.14}&
\|w_0(x)-w_0(t)\|\leq M_2 |x-t| \quad (M_2>0).
\end{align}
Clearly, $\wt \b(x)$ is continuous. Finally, formulas \eqref{2.6}, \eqref{2.9}, and \eqref{2.14} yield the inequality
\begin{align}\label{2.15}
\sup_{a\leq t<x\leq b}\| (x-t)^{-1}\wt \b(x)J\wt \b(t)^*\| \leq \wt M \quad  \big(\wt\b(x):=\b(x)w_0(x)\big).
\end{align}
Hence, in view of Lemma \ref{La}, we proved the following proposition.
\begin{Pn}\label{PnIneq} Let relations \eqref{2.5} and \eqref{2.6} hold, where $\b(x)$ is continuous on $[a,b]$.
Assume that parameter matrices $B$, $S(\xi)=S(\xi)^*$ and $\Pi(\xi)$ are given, that relations \eqref{1.3}, \eqref{1.17},
\eqref{1.14-}, and \eqref{2.0} hold, and  let $S(x)$ be invertible on $[a,b]$.
Then, the transformed Hamiltonian $\wt H(x)$ of the form \eqref{1.15} admits representation \eqref{2.10}, where $\wt \b$ 
is continuous and satisfies
the inequality \eqref{2.15}, which is similar to \eqref{2.6}. Moreover, the limit 
\begin{align}\label{2.7'}
\wt V(x,s)=\lim_{\eta \to +0}\big(\wt W(x, s+\I \eta)-\wt W(x, s-\I \eta)\big) \quad (a\leq s\leq x\leq b),
\end{align}
for $\wt W(x,z)$ of the form \eqref{2.3}, exists and its norm is uniformly bounded.
\end{Pn}
Lemma \ref{La} and further considerations in \cite{SaL4} show that $H(x)$ generates an interesting Riemann--Hilbert problem.
By virtue of Proposition \ref{PnIneq}, closely related Riemann--Hilbert problems are generated by the family of the
transformed Hamiltonians $\wt H(x)$.
\section{Characteristic functions \\ and spectral theory}\label{CharF}
\setcounter{equation}{0}
Let $\mfa\in \clb(\clh)$, $\clk\in \clb(\BC^m,\clh)$ and assume that the the identity
\begin{align} \label{c.1}&
\mfa-\mfa^*=\I \clk J \clk^*
\end{align}
holds. The   Liv\v{s}ic--Brodskii characteristic function \cite{BLiv, Liv} of the operator $\mfa$ is given by the formula
\begin{align} \label{c.2}&
\clw(z)=I_m-\I J\clk^*(\mfa - zI)^{-1}\clk .
\end{align}
Below, we consider the case of the operators $\mfa \in \clb\big(L^2_k(a,b)\big)$ of the form
\begin{align} \label{c.3}&
\big(\mfa  f\big)(x)=xf(x)+\I\int_a^x\b(x)J\b(t)^*f(t)dt \quad \big(\b(x)\in \BC^{k\times m} \big).
\end{align}
The operator $\mfa$ is a so called triangular Liv\v{s}ic model \cite{Liv, SaL5}
and its  characteristic function $W(b,z)$ is given by the formula \eqref{2.1}, where $H(x)=\b(x)^*\b(x)$ (see \cite{Liv, SaL4}).
In particular, taking into account that  $H(x)=\b(x)^*\b(x)$,  one easily derives:
\begin{align} \nn
(\mfa -z I)(z-x)^{-1}\b(x)W(x,z)=&-\b(x)W(x,z)
\\ \nn &
+\I \b(x)\int_a^xJ\b(t)^*\b(t)\frac{W(t,z)}{z-t}dt
\\ \nn
=& -\b(x)W(x,z)+\b(x)\int_a^x\frac{d}{dt}W(t,z)dt
\\ \nn
=& -\b(x).
\end{align}
That is, we have
\begin{align} \label{c.4}&
(\mfa -z I)^{-1}\b(x)=(x-z)^{-1}\b(x)W(x,z).
\end{align}
The equality $\clw(z)=W(b,z)$ for $\clw$ given by \eqref{c.2}, $\mfa$ of the form \eqref{c.3} and $\clk g=\b(x)g$ $(g\in \BC^m)$ easily follows from
\eqref{c.4}. 

Since $W(b,z)$ is the characteristic function of the operator $\mfa$ given by \eqref{c.3}, Theorem \ref{TmBDT} implies our next proposition.
\begin{Pn}\label{PnTChar} Let $H(x)=\b(x)^*\b(x)$ $\big(\b(x)\in\BC^{k\times m}\big)$ and assume that the equalities \eqref{2.0}
and the conditions of Theorem \ref{TmBDT} hold. Then, the characteristic function $\clw(z)$  of the operator $\mfa$ of the form
\eqref{c.3} and the characteristic function $\wt \clw(z)$  of the transformed operator $\wt \mfa$ given by
\begin{align} \label{c.5}&
\big(\wt \mfa  f\big)(x)=xf(x)+\I\int_a^x\wt \b(x)J\wt \b(t)^*f(t)dt \quad \big(\wt \b(x)=\b(x)w_0(x) \big).
\end{align}
are connected by the relation
\begin{equation} \label{c.6}
\wt \clw(z)=v(b, z) \clw(z)v(a, z)^{-1}.
\end{equation}
\end{Pn}

If the conditions of Lemma \ref{La} hold, then (according to \cite[Corollary 3.5]{SaL4}) the simple part of the operator  $\mfa$ of the form \eqref{c.3}
is linearly similar to a self-adjoint operator with an absolutely continuous spectrum. Hence, Proposition~\ref{PnIneq} yields the following result.
\begin {Cy}\label{SpT}
Let the conditions of Proposition \ref{PnIneq} be fulfilled. Then, the simple part of the operator  $\wt \mfa$ of the form \eqref{c.5}
is linearly similar to a self-adjoint operator with an absolutely continuous spectrum.
\end{Cy}
\begin{Rk}\label{Rkx} Let the operator $\mfa$ be given by \eqref{c.3}, where $J=I_m$, $k~=~m$ $($i.e., $\b(x)\in \BC^{m\times m})$,
$\b(x)\geq 0$ and the essential supremum of $\|\b(x)^2\|$ is bounded on $[a,b]$. Then, $\mfa$ is linearly similar  to the operator  $Q$
of the multiplication by $x$ in $\clh$, where $\clh$ is the closure of the linear manifold $\{g(x)=\b(x)f(x): f\in L^2_m(a,b)\}$ $($see \cite[Theorem 3.2]{SaL0}$)$.
If the conditions of Proposition \ref{PnTChar} are also fulfilled $($in addition to the conditions mentioned above in this remark$)$, 
we study transformed operators $\wt \mfa$ of the form \eqref{c.5},
where $\wt \b(x)=w_0(x)^*\b(x)w_0(x)$ and $w_0(x)w_0(x)^*\equiv I_m$ $($instead of the case $\wt \b(x)=\b(x)w_0(x)$  in the paper's
earlier  considerations$)$.
The operators $\wt \mfa$ $($with $\wt \b(x)=w_0(x)^*\b(x)w_0(x))$ are linearly similar to the
multiplication by $x$ in $\wt \clh$, where $\wt \clh$ is the closure of the linear manifold $\{\wt g(x)=\wt \b(x)f(x): f\in L^2_m(a,b)\}$  $($use again \cite[Theorem 3.2]{SaL0}$)$.
Moreover, in view of the equality $\wt \clw(z)=\wt W(b,z)$ and of Theorem \ref{TmBDT},  the characteristic functions $\wt \clw(z)\,$ of $\,\wt \mfa$
are connected with the characteristic function $\clw(z)\,$ of $\,\mfa$ via formula \eqref{c.6}.
\end{Rk}
\section{Explicit examples}\label{secEx}
\setcounter{equation}{0}
Let us  consider a simple case with a constant initial Hamiltonian $H$ of rank~1:
\begin{equation} \label{5.1}
\cli=[0,b], \,\, \xi=0,\,\, m=2, \,\, H(x) \equiv \b^* \b, \,\, \b \equiv [1 \quad \I], \quad  J=\left[\begin{array}{lr}0 & 1 \\ 1 & 0\end{array}
\right].
\end{equation}
It follows from \eqref{5.1} that $\b J\b^*=0$ and so formulas (\ref{1.1}) and (\ref{5.1}) yield
\begin{align} &\nn
\b W_x(x,z)=0, \quad \b J W_x(x,z)=2\I(z-x)^{-1} \b W(x,z), \quad
{\mathrm{i.e.,}}
\\ & \label{5.2}
\b W(x,z)= \b W(0,z)= \b, \quad \b J W(x,z)=-2 \I\big( \ln \,
(z-x) \big) \b +{\mathrm{const}},
\end{align}
where the normalisation condition  (\ref{1.2'}) was used (at $\xi=0$) and ``const" means some constant vector.
{\it Here and further, $z\not\in [0,b]$ and the branch of \, $\ln \,
(z-x)$ is fixed $($for the values of $x\in [0,b])$.}  Taking
into account (\ref{1.2'}) again, from the last  equality in
(\ref{5.2}), we derive
\begin{equation} \label{5.3}
\b J W(x,z)=2 \I \big( \ln \,
\frac{z}{z-x} \big) \b +\b J.
\end{equation}
According to \eqref{5.2} and \eqref{5.3}, we have
\begin{equation} \label{5.4}
W(x,z)=T^{-1}
\begin{bmatrix} \b \\ 2 \I \big( \ln \,
\frac{z}{z-x} \big) \b +\b J \end{bmatrix},
\end{equation}
where
\begin{equation} \label{5.4+}
T:=\begin{bmatrix}\b \\ \b J \end{bmatrix}=\begin{bmatrix}1 & \I \\ \I & 1 \end{bmatrix},
\quad TJT^*=2J, \quad T^{-1}=\frac{1}{2}\begin{bmatrix}1 & -\I \\-\I & 1 \end{bmatrix}.
\end{equation}
Clearly $W_{\pm}(x,s)$ do not exist for $s=0$ and $s=x$. According to \eqref{0.1} and \cite[(3.2)]{ALS07},
we have
\begin{equation} \label{5.5-}
W_+(x,s)=W_-(x,s)R(s)^2 \quad {\mathrm{for}}\quad R(s)=I_2+\pi J \b^*\b, \quad 0<s<x.
\end{equation}

In a similar to the construction of $W$ way, we recover $\Pi$. That is, from $\Pi_x=-\I A\Pi JH$ we obtain
\begin{equation} \label{5.5}
\Pi(x)J\b^*\equiv \Pi(0)J\b^*=g=\{g_i\}_{i=1}^n\in \BC^n.
\end{equation}
We also set $A=(B-x I_n)^{-1}$ and, using \eqref{5.5}, derive
\begin{equation} \label{5.6}
\frac{d}{dx}\Pi(x)  \b^* =-2 \I (B-xI_n)^{-1}\Pi(x) J \b^*=-2 \I (B-xI_n)^{-1}g.
\end{equation}
In the case
\begin{equation} \label{5.7}
A=(B-x I_n)^{-1}, \quad {\mathrm{where}}\quad  B= \mathrm{diag} \{ b_1, b_2, \ldots , b_n\}, \quad \s(B)\cap [0,b]=\emptyset,
\end{equation}
formula \eqref{5.6} implies that
\begin{equation} \label{5.8}
\Pi(x)  \b^* = 2\Big(\{ \I g_i \ln (b_i-x)\}_{i=1}^n +h \Big) \qquad
(h \in \BC^n) ,
\end{equation}
where the branch of $\ln(b_i-x)$ $(x\in[0,b])$ is fixed for each $i$ $(1\leq i \leq n)$. It follows from \eqref{5.4+}, \eqref{5.5}, and \eqref{5.8} that
\begin{equation} \label{5.9}
\Pi(x)   =\frac{1}{2}\begin{bmatrix} 2\Big(\{ \I g_i \ln (b_i-x)\}_{i=1}^n +h \Big) & g\end{bmatrix}
T,
\end{equation}
where $h$ and $g$ are easily expressed in terms of the parameter matrix $\Pi(0)$. That is, we have:
\begin{equation} \label{5.10}
g=\Pi(0)   \begin{bmatrix}- \I \\ \ 1 \end{bmatrix}, \quad h=\frac{1}{2}\Pi(0)\begin{bmatrix}1 \\ \ -\I \end{bmatrix}-\{ \I g_i \ln (b_i)\}_{i=1}^n .
\end{equation}
Since  $b_k \not\in [0, b] $, $\Pi(x)$ is
well defined on $[0, b]$. Taking into account (\ref{5.4+}) and \eqref{5.9} we obtain
\begin{align} &\Pi(x)J \Pi(x)^*
 \label{5.11}
=\Big(\{ \I g_i \ln (b_i-x)\}_{i=1}^n +h \Big)g^*+g\Big(\{ \I g_i
\ln (b_i-x)\}_{i=1}^n +h \Big)^*.
\end{align}
Using \eqref{5.7}, the matrix function $S=\{S_{ij}\}_{i,j=1}^n$ may be derived from the identity \eqref{1.6}:
\begin{equation} \label{5.12}
S_{ij}(x)=\frac{(b_i-x)(\ov{b_j}-x)}{\I (b_i-\ov{b_j})}\Big(\Pi (x)J \Pi(x)^*\Big)_{ij} \quad {\mathrm{for}} \quad b_i\not=\ov{b_j},
\end{equation}
where $\Big(\Pi (x)J \Pi(x)^*\Big)_{ij} $ are the entries of $\Pi (x)J \Pi(x)^*$.
Now, in view of (\ref{1.11}) and (\ref{5.5}) we obtain
\begin{equation}
\label{5.13-}
\b w_0(x)= \b - \I g^*S(x)^{-1}(B-x I_n) \Pi(x),
\end{equation}
which, taking into account (\ref{1.15}) and $H(x) \equiv \b^* \b$, implies
\begin{equation} \label{5.13}
\wt H(x)=  \big( \b - \I g^*S(x)^{-1}(B-x I_n) \Pi(x)
\big)^*
\big( \b - \I g^*S(x)^{-1}(B-x I_n) \Pi(x) \big)
\end{equation}
for the GBDT-transformed Hamiltonian $\wt H(x)$.
From \eqref{1.11} and \eqref{5.8}, it follows that
\begin{equation} \label{5.14}
\b Jw_0(x)=\b J- 2 \I \Big(\{ \I g_i \ln (b_i-x)\}_{i=1}^n +h \Big)^*S(x)^{-1}(B-x I_n) \Pi(x).
\end{equation}
Relations \eqref{5.13-} and \eqref{5.14} imply that
\begin{equation} \label{5.15}
w_0(x)=I_2-\I T^{-1}\begin{bmatrix}g^* \\ 2  \Big(\{ \I g_i \ln (b_i-x)\}_{i=1}^n +h \Big)^*\end{bmatrix}S(x)^{-1}(B-x I_n) \Pi(x).
\end{equation}
\begin{Rk}\label{RkEx}
Relations \eqref{1.16} and \eqref{1.7}, \eqref{1.10}, \eqref{1.17+} together with the explicit formulas 
\eqref{5.9}--\eqref{5.12} and \eqref{5.15} give explicit expressions of $\wt W(x,z)$ for the case of
$\wt H(x)$ given by \eqref{5.13} $($and $\xi=a=0)$. Since $\b J\b^*=0$, we see that \eqref{2.5} and \eqref{2.6}
holds. Thus, our example satisfies conditions of  Proposition \ref{PnIneq}.
Moreover, according to Proposition \ref{PnTChar} the constructed matrix function $\wt W(b,z)$ is the
characteristic matrix function of the operator $\wt \mfa$ given by \eqref{c.5}, where $a=0$ and $\wt \b(x)$
is given by \eqref{5.13-}.
\end{Rk}
Let us consider the simplest case $n=1$ in a more detailed way.
\begin{Ee}
Let $n=1$. Then, $g$, $h$, and $B$ are scalars  and we assume that $B
\not\in \BR$. Rewrite (\ref{5.11}) as
\begin{equation} \label{e1}
\Pi(x)J \Pi(x)^*=\I |g|^2 \big( \ln (B-x)- \ov{\ln (B-x)} \big)+ h
\ov{g}+g \ov{h}.
\end{equation}
Note that, according to  \eqref{5.10}, there is a simple correspondence
between the pairs $g$ and $h$ and the parameter matrix  $\Pi(0)$ so that one may fix
$g$ and $h$ instead of $\Pi(0)$. Next, in
view of (\ref{5.12}) and \eqref{e1}, we obtain
\begin{equation} \label{e2}
S(x)= \Big( |g|^2 \big( \ln (B-x)- \ov{\ln (B-x)} \big) -\I h
\ov{g}-\I g \ov{h} \Big) \frac{(B-x)(\ov{B}-x)}{B- \ov{B}}.
\end{equation}
Since $B \not\in \BR$, the condition $S(x)\not=0$ for any $x\in [a,b]$
$($the condition of invertibility of $S(x))$ may be written down as
\begin{equation} \label{e3}
g\not=0, \quad \Im  \big( \ln (B-x)\big)\not=\Re(h/{g}).
\end{equation}
Under conditions \eqref{e3}, taking into account \eqref{5.9}, \eqref{5.13-}, and \eqref{e2}, we derive
for $\wt \b(x)= \b w_0(x)$ the equality
\begin{align} \nn
\wt \b(x)=&\b-\frac{\I \ov{g}(B- \ov{B})}{2(\ov{B}-x)\Big( |g|^2 \big( \ln (B-x)- \ov{\ln (B-x)} \big) -\I h
\ov{g}-\I g \ov{h} \Big)}
\\ & \label{e4}
\times \begin{bmatrix} 2\Big( \I g\ln (B-x) +h \Big) & g\end{bmatrix}T.
\end{align}
Clearly, $\Im  \big( \ln (B-x)\big)\not=0$. Hence, \eqref{e3} holds for $g\not=0$ and $h=0$. In this case, formula \eqref{e4} is simplified
and takes the form
\begin{align}  & \label{e5}
\wt \b(x)=\b-\frac{(B- \ov{B})}{4(\ov{B}-x)\Im  \big( \ln (B-x)\big)}
 \begin{bmatrix} 2 \I \ln (B-x) & 1\end{bmatrix}T.
\end{align}
Using \eqref{e4} or \eqref{e5}, one obtains the GBDT-transformed Hamiltonian
$\wt H(x)=\wt \b(x)^*\wt \b(x)$ and the GBDT-transformed triangular model
\begin{align} \label{e6}&
\big(\wt \mfa  f\big)(x)=xf(x)+\I\int_a^x\wt \b(x)J\wt \b(t)^*f(t)dt .
\end{align}
In view of \eqref{5.15} $($using also  the formulas \eqref{5.9} and \eqref{e2}$)$, we have
\begin{align}\nn
w_0(x)=& I_2-\frac{\I (B- \ov{B})}{2(\ov{B}-x)\Big( |g|^2 \big( \ln (B-x)- \ov{\ln (B-x)} \big) -\I h
\ov{g}-\I g \ov{h} \Big)}
\\ & \label{e7}
\times  T^{-1}\begin{bmatrix}g^* \\ 2  \big( \I g \ln (B-x) +h \big)^*\end{bmatrix}\begin{bmatrix} 2\Big( \I g\ln (B-x) +h \Big) & g\end{bmatrix}T.
\end{align}
In a similar way, formulas \eqref{1.7} and \eqref{1.17+} yield
\begin{align}\nn
w_A(x,z)=& I_2-\frac{\I(x-z) (B- \ov{B})}{2(B-z)(\ov{B}-x)\Big( |g|^2 \big( \ln (B-x)- \ov{\ln (B-x)} \big) -\I h
\ov{g}-\I g \ov{h} \Big)}
\\ & \label{e8}
\times  T^{-1}\begin{bmatrix}g^* \\ 2  \big( \I g \ln (B-x) +h \big)^*\end{bmatrix}\begin{bmatrix} 2\Big( \I g\ln (B-x) +h \Big) & g\end{bmatrix}T.
\end{align}
It follows from  \eqref{1.10},  \eqref{e7},  \eqref{e8}, and the equality $w_0(x)^{-1}= Jw_0(x)^*J$ that
\begin{align} \label{e9}&
v(x,z)=\frac{1}{B-z}\Big((x-z)I_2+(B-x)Jw_0(x)^*J\Big).
\end{align}
Recall that $a=0$ in this section and that the characteristic function $\clw(z)$ of the operator $\mfa$ equals $W(b,z)$.   Thus, finally,
relations \eqref{c.6}, \eqref{5.4}, and \eqref{e9} give us an explicit formula the  characteristic function $\wt \clw(z)$
 of the GBDT-transformed operator $\wt \mfa$.
\end{Ee}

\end{document}